\documentclass[reqno,12pt]{amsart}
\usepackage{a4wide}
\usepackage{amsmath}
\usepackage{amssymb}
\usepackage{amsthm}
\usepackage[frenchb]{babel}
\usepackage{mathrsfs}

\newtheorem{theo}{Th\'eor\`eme}

\newtheorem{coro}{Corollaire}

\newtheorem{lem}{Lemme}

\newtheorem*{obs}{Observation 2}

\begin{document}

\title[]{D\'emonstration de l'Observation 2 d'Almkvist et Zudilin}
\author[]{C. Krattenthaler et T. Rivoal}
\date{\today}
\maketitle

\markboth{}{}

Le but de cette note est de donner une  preuve de l'Observation 2 
d'Almkvist et Zudilin~\cite[p. 487]{az}. Il 
s'agit de l'observation num\'erique suivante. Posons
$$
A_n= \sum_{0\le j,k \le n}\binom{n}{j}^2\binom{n}{k}^2\binom{n+j}{n}\binom{n+k}{n}\binom{j+k}{n}
$$
et 
$$
w_0(z)=\sum_{n=0}^{\infty} A_n z^n = 1+12z+804z^2+88680z^3+\cdots \in \mathbb{Z}[[z]].
$$
\begin{obs} Il existe $u_0(z) \in \mathbb{Z}[[z]]$ telle que $w_0(z)=u_0(z)^2.$
\end{obs}

On va utiliser le crit\`ere suivant de Heninger et al.~\cite{sloane}. On note $\mathcal{P}_n$ l'ensemble de s\'eries 
formelles $F(z)\in 1+z\mathbb{Z}[[z]]$ telles que $F(z)^{1/n}\in 1+z\mathbb{Z}[[z]]$ pour un entier $n\ge 1.$

\begin{lem} Posons $\mu_n=n\prod_{p\vert n} p$. On a  
$$
F(z) \in \mathcal{P}_n \Longleftrightarrow F(z) \;(\textup{mod} \,\mu_n) \in \mathcal{P}_n.
$$
\end{lem}

Dans le cas de la s\'erie $w_0(z)$, on a $n=2$ et $\mu_2=4.$ Il suffit donc de prouver que 
$w_0(z) \;(\textup{mod} \,4) \in \mathcal{P}_2$. Or nous allons prouver que 
$w_0(z)\;(\textup{mod} \,4) = 1$, qui est bien dans $\mathcal{P}_2$.

\begin{theo} 
\label{theo:1}
Pour des entiers $j,k,n\ge 0$, posons
$$
a_n(j,k)=\binom{n}{j}^2\binom{n}{k}^2\binom{n+j}{n}\binom{n+k}{n}\binom{j+k}{n}.
$$ 

\noindent Pour tous entiers $n\ge 1$ et $j,k\ge 0$, on a

$(a)$ 
$v_2\big(a_n(j,k)\big)=1$ 
si, et seulement si, $\{j,k\}=\{0,n\}$ et $n$ est une puissance de $2$.

$(b)$ 
$v_2\big(a_n(j,k)\big)\ge 2$
sinon.
\end{theo}

Pour $n\ge 1$, on a 
$$
A_n= a_n(0,n)+a_n(n,0)+
{\sum_{0\le j,k \le n}}\kern-5pt{}^{\displaystyle\prime} \kern5pt
a_n(j,k), 
$$
o\`u le prime signifie que la sommation exclut les couples $(j,k)=(0,n), (n,0)$.
Comme $a_n(0,n)=\binom{2n}{n}=a_{n}(n,0)$ est toujours pair (le pire cas \'etant donn\'e par (a)), on en d\'eduit 
le corollaire suivant.

\begin{coro} Pour tout entier $n\ge 1$, on a $v_2\big(A_n\big)\ge 2$. En particulier, 
$w_0(z)\;(\textup{mod} \,4) = 1$ et l'Observation 2 est vraie.
\end{coro}

\begin{proof}[D\'emonstration du Th\'eor\`eme~\ref{theo:1}] 
Soit $n\ge1$.
Il n'y a rien \`a prouver si $j+k<n$ et on suppose maintenant que $j+k\ge n.$
\medskip

{\bf Premier cas : $\{j,k\}=\{0, n\}$.} On a alors 
$a_n(0,n)=a_n(n,0)=\binom{2n}{n}$. Or 
$$
v_2\bigg(\binom{2n}{n}\bigg) = 
\sum_{\ell=1}^{\infty} \bigg(\bigg \lfloor \frac{2n}{2^\ell} \bigg\rfloor - 
2\bigg \lfloor \frac{n}{2^\ell} \bigg\rfloor\bigg) 
= n- \sum_{\ell=1}^{\infty} \bigg \lfloor \frac{n}{2^\ell} \bigg\rfloor.
$$
Il est facile de v\'erifier que l'expression \`a droite est $\ge 2$ sauf si $n$ est une puissance de $2$, auquel cas elle vaut $1$.

\medskip

{\bf Deuxi\`eme cas : $\{j,k\} \neq \{0, n\}$.} La condition $j+k\ge n\ge 1$ montre que dans  
ce cas on a alors forc\'ement $j\ge 1 $ et $k\ge 1.$ 

On va tout d'abord montrer que pour entiers $m,p$ tels que $1\le p\le m$, on a 
\begin{equation}\label{eq:1}
v_2\bigg(\binom{m}{p}\binom{m+p}{m}\bigg) \ge 1.
\end{equation}
En effet, on a 
\begin{align*}
v_2\bigg(\binom{m}{p}\binom{m+p}{m}\bigg) 
&= \sum_{\ell=1}^{\infty}\bigg( \bigg \lfloor \frac{m+p }{2^\ell} \bigg\rfloor - 
2\bigg \lfloor \frac{p}{2^\ell} \bigg\rfloor - \bigg \lfloor \frac{m-p}{2^\ell} \bigg\rfloor\bigg)
\\
&\ge \bigg \lfloor \frac{m+p }{2^L} \bigg\rfloor - 
2\bigg \lfloor \frac{p}{2^L} \bigg\rfloor - \bigg \lfloor \frac{m-p}{2^L} \bigg\rfloor 
= \bigg \lfloor \frac{m+p }{2^L} \bigg\rfloor  - \bigg \lfloor \frac{m-p}{2^L} \bigg\rfloor \ge 0,
\end{align*}
o\`u $L$ est l'entier $\ge 1$ tel que $2^{L-1}\le p<2^L.$
Supposons alors que $\lfloor \frac{m+p }{2^L} \rfloor  - \lfloor \frac{m-p}{2^L} \rfloor = 0$. On en d\'eduit que 
$$
0\le \frac{m+p }{2^L} - \frac{m-p}{2^L}<1
$$
et donc que $p<2^{L-1}$, contrairement \`a l'hypoth\`ese. Donc  $\lfloor\frac{m+p }{2^L} \rfloor  
- \lfloor \frac{m-p}{2^L} \rfloor \ge 1$, ce qui prouve~\eqref{eq:1}.

On applique maintenant~\eqref{eq:1} aux deux cas $(m,p)=(n,j)$ et  $(m,p)=(n,k)$ pour obtenir que 
$$
v_2\bigg(\binom{n}{j}\binom{n}{k}\binom{n+j}{n}\binom{n+k}{n}\bigg) \ge 2.
$$
{\em A fortiori}, on a donc 
$v_2\big(a_n(k,j)\big)\ge 2$ lorsque $\{j,k\}\neq \{0,n\}.$

\medskip

Ceci termine la preuve du th\'eor\`eme.
\end{proof}

\def\refname{Bibliographie}

\end{document}